\documentclass[12pt,a4paper]{article}
\usepackage{amssymb,enumerate,verbatim}
\usepackage{amsgen,amsmath,amsthm,amstext,amsbsy,amsopn,amsfonts,amssymb}

\newtheorem{Def}{Definition}
\newtheorem{Thm}{Theorem}
\newtheorem{Prop}{Proposition}
\newtheorem{Lem}{Lemma}
\newtheorem{Cor}{Corollary}
\newtheorem{Ex}{Example}
\newtheorem{Rem}{Remark}
\newcommand\Z{\mathbb{Z}}

\newcommand\Q{\mathbb{Q}}
\newcommand\R{\mathbb{R}}

\newcommand\F{\mathbb{F}}
\newcommand\N{\mathbb{N}}
\DeclareMathOperator{\Ima}{Im}

\title{Expanders and box spaces}
\author{Ana KHUKHRO and Alain VALETTE}

\begin{document}
\maketitle

\begin{abstract} We consider box spaces of finitely generated, residually finite groups $G$,  and try to distinguish them up to coarse equivalence. We show that, for $n\geq 2$, the group $SL_n(\Z)$ has a continuum of box spaces which are pairwise non-coarsely equivalent expanders. Moreover, varying the integer $n\geq 3$, expanders given as box spaces of $SL_n(\Z)$ are pairwise inequivalent; similarly, varying the prime $p$, expanders given as box spaces of $SL_2(\Z[\sqrt{p}])$ are pairwise inequivalent.

A strong form of non-expansion for a box space is the existence of $\alpha\in]0,1]$ such that the diameter of each component $X_n$ satisfies $diam(X_n)=\Omega(|X_n|^\alpha)$. By \cite{BrTo}, the existence of such a box space implies that $G$ virtually maps onto $\Z$: we establish the converse. For the lamplighter group $(\Z/2\Z)\wr\Z$ and for a semi-direct product $\Z^2\rtimes\Z$, such box spaces are explicitly constructed using specific congruence subgroups.

We finally introduce the full box space of $G$, i.e. the coarse disjoint union of all finite quotients of $G$. We prove that the full box space of a group mapping onto the free group $\F_2$ is not coarsely equivalent to the full box space of an $S$-arithmetic group satisfying the Congruence Subgroup Property. 
\end{abstract}

\section{Introduction}

If $(X,d_X),(Y,d_Y)$ are metric spaces, a map $f:X\rightarrow Y$ is a {\bf coarse embedding} if for any two sequences $(x_n)_{n>0}, (y_n)_{n>0}$ in $X$: 
$$d_X(x_n,y_n)\rightarrow +\infty \Longleftrightarrow d_Y(f(x_n),f(y_n))\rightarrow +\infty.$$
Moreover, $f$ is a {\bf coarse equivalence} if $f$ is almost surjective, i.e. there exists $R\geq 0$ such that $Y$ is the $R$-neighbourhood of $f(X)$.

Let $X=\coprod_{n>0} X_n$ be a disjoint union of finite, connected, $d$-regular graphs. The $X_n$'s are called {\bf components} of $X$. Endow $X$ with a metric $d$ (well-defined up to coarse equivalence) inducing the original metric on each $X_n$ and such that $d(X_m,X_n)\geq \max\{diam(X_m),diam(X_n)\}$ for $m\neq n$: we call $X$ the {\bf coarse disjoint union} of the $X_n$'s. We say that $X$ is an {\bf expander} if, for some $\epsilon>0$, we have $h(X_n)\geq\epsilon$ for every $n>0$, where $h(X_n)$ denotes the Cheeger constant (or isoperimetric constant) of the graph $X_n$.

The question of distinguishing expanders up to coarse equivalence was first considered by M. Mendel and A. Naor (\cite{MeNa}, Theorem 9.1; see also Ostrovskii \cite{Ost}, Theorem 5.76); they constructed an expander $X$ with unbounded girth (i.e. $g(X_n)\rightarrow +\infty$), an expander $Y$ containing many short cycles, and proved that $X$ does not coarsely embed into $Y$. This was dramatically generalized recently by D. Hume \cite{Hume}, who constructed a continuum of expanders, all with unbounded girth, and proved that pairwise they do not embed coarsely into one another. This leaves the question of distinguishing expanders with bounded girth up to coarse equivalence.

We attack this question using box spaces of residually finite groups, of which we now recall the definition. Let $G$ denote a finitely generated, residually finite group. A decreasing sequence $(N_k)_{k>0}$ of finite index, normal subgroups of $G$ with trivial intersection will be called a {\bf filtration} of $G$.

\begin{Def}If $(N_k)_{k>0}$ is a filtration of $G$, the corresponding box space $\square_{(N_k)}G$ is the coarse disjoint union of the $G/N_k$'s. 
\end{Def}

Fix a finite generating set $S$. Each $G/N$ is endowed with the word metric coming from $S$ (i.e. is viewed as Cayley graph with respect to $S$). We then endow $\square_{(N_k)}G$ with a metric as above. It follows from Proposition 2 in \cite{Ana2012} that, up to coarse equivalence, $\square_{(N_k)}G$ does not depend on the choice of the finite generating set $S$. 



Box spaces are known to be strongly connected to expanders since Margulis \cite{Mar75} used box spaces of property (T) groups to give the first explicit construction of expanders. We will prove: 

\begin{Thm}\label{SLm(pk)} Let $m\geq 2$. Denote by $\Gamma(N)=\ker[SL_m(\Z)\rightarrow SL_m(\Z/N\Z)]$ the congruence subgroup. For $s\geq 1$ a real number, set $N_k(s)=2^{[ks]}$ (with $[.]$ denoting the floor function). The box spaces $\square_{(\Gamma(N_k(s)))} SL_m(\Z)$ are expanders which are pairwise not coarsely equivalent. In particular (taking $m=3$), there exists a continuum of expanders with geometric property (T) of \cite{WiYu}.
\end{Thm}

\begin{Thm}\label{Hilbert} For prime $p$, let $X_p$ be any box space of the Hilbert modular group $SL_2(\Z[\sqrt{p}])$. For $n\geq 3$, let $Y_n$ be any box space of $SL_n(\Z)$.
\begin{enumerate}
\item $X_p$ is an expander and, for distinct primes $p, q$, there is no coarse equivalence between $X_p$ and $X_q$.
\item For $n\geq 3$, $Y_n$ is an expander and, for $m\neq n$, there is no coarse equivalence between $Y_m$ and $Y_n$.
\item For $n\geq 3$ and $p$ prime, there is no coarse embedding of $Y_n$ into $X_p$.
\end{enumerate}
\end{Thm}

It is known (see Proposition 7.3.11 in \cite{Lub}) that, if $X=\coprod_n X_n$ is an expander, then the diameter of $X_n$ is at most logarithmic in the number of vertices of $X_n$. Therefore, requiring the diameter of $X_n$ to grow like a small power of $|X_n|$ is a strong form of non-expansion:

\begin{Def} Fix $\alpha\in ]0,1]$. We say that $\square_{(M_k)} G$ has {\bf property $D_\alpha$} if there exists $K>0$ (only depending on $S$) such that, for every component $G/M$ of $\square_{(M_k)} G$:
$$diam(G/M)\geq K.|G/M|^\alpha.$$
\end{Def}

We prove that property $D_\alpha$ is a coarse invariant of box spaces. Moreover, in the extreme case $\alpha=1$, we have the following result.

\begin{Thm}\label{virtuallycyclic} The following are equivalent:
\begin{enumerate}
\item[i)] $\square_{(M_k)} G$ satisfies property $D_1$;
\item[ii)] $G$ is virtually cyclic.
\end{enumerate}
\end{Thm}

Property $D_\alpha$ is related to a property recently introduced by E. Breuillard and M. Tointon \cite{BrTo}: if $H$ is a finite group generated by some subset $S$, say that the corresponding Cayley graph is {\bf $\alpha$-almost flat} if 
$$diam_S(H)\geq (\frac{|H|}{|S|})^\alpha.$$
Decreasing $\alpha$ if necessary, and discarding finitely many $N_k$'s, we see that property $D_\alpha$ for $\square_{(N_k)}G$ is equivalent to all $G/N_k$'s being $\alpha'$-almost flat for some $\alpha'\leq\alpha$. It is known that box spaces of virtually nilpotent groups have $D_\alpha$ (\cite{BrTo}, Proposition 4.7 and Remark 4.3). An important structural consequence of property $D_\alpha$ was obtained by E. Breuillard and M. Tointon (combine Corollary 1.7 with Lemma 5.1 in \cite{BrTo}): if $\square_{(N_k)}G$ has $D_\alpha$, then $G$ virtually maps onto $\Z$. Note that the converse is false: the free group $\F_2$ maps onto $\Z$, but admits box spaces which are expanders. However the converse holds if one just requests SOME box space with property $D_\alpha$:

\begin{Thm}\label{Breuillard} Let $G$ be a finitely generated, residually finite group. TFAE:
\begin{enumerate}
\item[a)] $G$ virtually maps onto $\Z$;
\item[b)] for some $\alpha>0$ and some filtration $(N_k)_{k>0}$ of $G$, the box space $\square_{(N_k)} G$ has property $D_\alpha$.
\end{enumerate}
\end{Thm}

For the lamplighter group $(\Z/2\Z)\wr\Z$ and for the semi-direct product $\Z^2\rtimes\Z$ (where $\Z$ acts by powers of $A=\left(\begin{array}{cc}1 & 1 \\1 & 0\end{array}\right)$), we explicitly construct such box spaces using specific congruence subgroups.

Box spaces of a given group $G$ depend rather crucially on the choice of a filtration in $G$; for example, the free group has some box spaces which are expanders, and some others which are not. To get something canonically associated with $G$, we define the {\bf full box space} $\square G$ as the disjoint union of all finite quotients of $G$:
$$\square G=\coprod_{N\triangleleft G,\;[G:N]<+\infty} G/N,$$
which we endow with a metric (well-defined up to coarse equivalence) coming from some finite generating of $G$, exactly as in the previous case.

\begin{Thm}\label{fullvirtuallycyclic} The full box space of a virtually cyclic, infinite group is coarsely equivalent to the full box space of $\Z$, i.e. the disjoint union of all cycles $C_n$.
\end{Thm}

We show by example that the analogous statement is no longer true for groups which are virtually $\Z^2$. By relating the coarse geometry of $\square G$ with the asymptotics of the counting function of normal subgroups of a given in $G$, we give a sufficient condition for distinguishing full box spaces:

\begin{Thm}\label{CSP} Suppose that $G$ maps onto the free group $\mathbb{F}_2$, and that $H=\mathbb{G}(O_S)$, an $S$-arithmetic subgroup in a simple, classical\footnote{More precisely: $\mathbb{G}$ is NOT of type $G_2, F_4$ or $E_8$.} algebraic group $\mathbb{G}$ defined over some number field $k$. If $H$ has the congruence subgroup property (CSP), then $\square G$ and $\square H$ are NOT coarsely equivalent.
\end{Thm}

The paper is organized as follows: in Section 2 we prove the main lemma, i.e. that coarse equivalences between box spaces are actually families of quasi-isometries with uniform constants between components. Theorem \ref{Hilbert} is proved in Section 3, while property $D_\alpha$ is dealt with in Section 4. In Section 5 we discuss full box spaces and prove Theorem \ref{fullvirtuallycyclic}. Section 6 is devoted to examples: apart from proving Theorem \ref{CSP}, we provide examples of quasi-isometric (even commensurable!) groups having non-coarsely equivalent full box spaces; and of non-commensurable groups admitting coarsely equivalent box spaces. We finish with some open problems in Section 7.

\medskip
{\bf Acknowledgements:} We thank V. Alekseev, K. Bou-Rabee, E. Breuillard, M. Finn-Sell, A. Lubotzky, N. Nikolov, and R. Tessera for useful input at various stages of the project.

\section{The main lemma}

We will use the standard fact that, for two metric spaces $(X,d_X),(Y,d_Y)$, a map $f: X\rightarrow Y$ is a coarse embedding if and only if there exist {\bf control functions} $\rho_\pm:\R^+\rightarrow\R^+$, with $\lim_{t\rightarrow+\infty}\rho_-(t)=+\infty$, and 
$$\rho_-(d_X(x,y))\leq d_Y(f(x),f(y))\leq\rho_+(d_X(x,y))$$ for every $x,y\in X$.


Let $X= \coprod_{k>0} X_k$ and $Y= \coprod_{k>0} Y_k$ be coarse disjoint unions of graphs of strictly increasing diameter; let $f: X \rightarrow Y$ be a coarse equivalence with control functions $\rho_\pm$, with $\Ima\;f$ being $C$-dense in $Y$. By increasing $C$ if necessary, we may assume $C\geq \rho_+(1)$. 

Set $t_C=\inf\{t>0:\rho_-(t)>C\}$; we define finite subsets $F_X\subset X$ and $F_Y\subset Y$ by 
$$F_X=\coprod_{k: \,d_{X}(X_k,(X_k)^c)\leq t_C} X_k;\;F_Y=\coprod_{k:\,d_{Y}(Y_k,(Y_k)^c)\leq C} Y_k.$$

\begin{Lem}\label{main} $f$ builds up a bijection between a co-finite set of components of $X$ and a co-finite set of components of $Y$. Moreover, there exists a constant $A>0$ such that, for every component $X_k$ of $X$, the restriction $f|_{X_k}$ satisfies: for every $x,x'\in X_k$,
$$\frac{1}{A}d_{X}(x,x')-A\leq d_{Y}(f(x),f(x'))\leq A\cdot d_{X}(x,x').$$
In particular $f|_{X_k}$ is a $(A,A)$-quasi-isometry.
\end{Lem}

{\bf Proof:} \begin{itemize}
\item Take $x,y\in X_k$, at distance 1, and assume that $f(x)$ belongs to some component $Y_\ell \subset Y-F_Y$. Then $f(y)$ also belongs to $Y_\ell$ as $d_{Y}(f(x),f(y))\leq \rho_+(1)= C$.
\item For $Y_\ell\subset Y-F_Y$, there exists $X_k$ such that $f(X_k)$ meets $Y_\ell$. By the previous step, and using connectedness of $X_k$, we have $f(X_k)\subset Y_\ell$. Observe that $Y_\ell$ is the $C$-neighbourhood of $f(X_k)$.
\item If $X_k\subset X-F_X$ and $f(X_k)\subset Y_\ell\subset Y - F_Y$, then $X_k$ is the only component of $X$ mapping to $Y_\ell$. Indeed, assuming $f(X_j)\subset Y_\ell$, we get $C\geq d_{Y}(f(X_j),f(X_k))\geq\rho_-(d_{X}(X_j,X_k))$, i.e. $d_{X}(X_j,X_k)\leq t_C$, hence $X_j=X_k$.
\item Consider those components of $X-F_X$ that map to $Y- F_Y$; by the previous step, they map injectively. They are co-finite: since $f$ has finite fibers, there are finitely many components of $X$ mapping to $F_Y$. The image of the components of $X-F_X$ that map to $Y- F_Y$, can also be described as the set of components of $Y-F_Y$ that meet $f(X - F_X)$. This set is co-finite: indeed by the second step of the proof, a component of $Y-F_Y$ disjoint from $f(X - F_X)$, is attained by some component of $X$ necessarily in $F_X$, and there are finitely many such components.

\item Finally, to prove the metric statement: for every $x,x'\in X_k$, we have $d_{Y}(f(x),f(x'))\leq \rho_+(1)d_{X}(x,x')$. To get the lower bound, let $g:Y\rightarrow X$ be a coarse inverse for $f$, with control functions $\rho'_-,\rho'_+$. Then (as above): $d_{X}(g(f(x)),g(f(x')))\leq \rho'_+(1).d_{Y}(f(x),f(x'))$. Let $K>0$ be such that every $x\in X_k$ is $K$-close to $g(f(x))$. Then $d_{X}(x,x')-2K\leq d_{X}(g(f(x)),g(f(x')))$, and the result follows.
\end{itemize}\hfill$\square$

\begin{Ex} If $f$ is only assumed to be a coarse embedding, it is not true that $f$ is a quasi-isometric embedding. To see it, set on the one hand $G=\Z$ and $M_k=2^k\Z$; on the other hand set $H=Heis(\Z)$, the $3\times 3$ Heisenberg group over $\Z$, and $N_k=\ker(Heis(\Z)\rightarrow Heis(\Z/2^k\Z))$, the kernel of reduction modulo $2^k$. Let $f:\square_{(M_k)} G\rightarrow \square_{(N_k)} H$ be the map obtained by viewing $\Z/2^k\Z$ as the center of $Heis(\Z/2^k\Z)$. Because $\Z$, viewed as the center of $H$, is polynomially distorted in $H$, the map $f$ is a coarse embedding but not a quasi-isometric embedding. Specifically, for the element $2^{k-1}\in\Z/2^k\Z$, the word length of $f(2^{k-1})\in Heis(\Z/2^k\Z)$ is of the order of $2^{\frac{k-1}{2}}$.
\end{Ex}

Lemma \ref{main} motivates the following definition:

\begin{Def} An almost permutation of $\N$ is a bijection between two co-finite subsets of $\N$. If $\alpha$ is an almost permutation, we say that $\alpha$ has bounded displacement if there exists $N$ such that $|\alpha(k)-k|\leq N$ for $k\gg 0$.
\end{Def}

We now make Lemma \ref{main} more precise in the case of box spaces. Note that in what follows, we will assume $0\notin \mathbb{N}$. We will need a lemma.

\begin{Lem}\label{permut} Let $\alpha$ be a permutation of $\N$. Assume that there exists $N\geq 1$ such that:
$$\forall k,\ell\in\N: \ell\geq k+N \Rightarrow \alpha(\ell)>\alpha(k).$$
Then $\alpha$ has bounded displacement.
\end{Lem}

{\bf Proof:} We first prove that $\alpha(k)\geq k-N$ for every $k\geq N$. For $1\leq i\leq k-N$, we have by assumption $\alpha(i)<\alpha(k)$. Since there are $k-N$ possibilities for $i$, this forces $\alpha(k)\geq k-N$.

We now prove that $\alpha(k)\leq k+N$ for every $k$. Assume by contradiction that this is false, and let $k$ be the smallest integer with $\alpha(k)>k+N$. Then $\alpha(i)\leq i+N$ for $1\leq i\leq k-1$. Set $F=:\{1,2,...,k+N\}\backslash \{\alpha(1),...,\alpha(k-1)\}$. 
For each $\ell$ in $\alpha^{-1}(F)$, we have $\ell>k-1$. But since $|\alpha^{-1}(F)|=|F|=N+1$, there must exist some $\ell\in \alpha^{-1}(F)$ with $\ell\geq k+N$, and $\alpha(\ell)\leq k+N<\alpha(k)$, in contradiction with the assumption.
\hfill$\square$

\medskip
Let $(N_k)_{k>0}$ be a filtration of the residually finite group $G$. We say that it is a {\it strict filtration} if the sequence $(N_k)_{k>0}$ is strictly decreasing.

\begin{Prop}\label{bdddispl} Let $(M_k)_{k>0}, (N_k)_{k>0}$ be strict filtrations of the residually finite groups $G,H$ respectively. Let $f: \square_{(M_k)} G\rightarrow \square_{(N_k)}H$ be a coarse equivalence and, as in Lemma \ref{main}, let $\alpha$ be the bijection induced between two co-finite sets of components, viewed as an almost permutation of $\N$. Then $\alpha$ has bounded displacement.
\end{Prop}

{\bf Proof:} Discarding finitely many initial $M_k$'s and $N_k$'s and renumbering accordingly, we may assume that, for every $k>0$, the map $f$ maps $G/M_k$ to $H/N_{\alpha(k)}$, and we view $\alpha$ as a permutation of $\N$. Let $B_H(R)$ be the ball of radius $R$ at the identity of $H$. Since $H/N_{\alpha(k)}$ is a $C$-neighbourhood of $f(G/M_k)$, we have
\begin{equation}\label{volume}
|H/N_{\alpha(k)}|\leq |B_H(C)|\cdot |f(G/M_k)|\leq |B_H(C)|\cdot |G/M_k|
\end{equation}
so, by symmetry, there exists a constant $K\geq 1$ such that, for all $k\in\N$:
\begin{equation}\label{1symm}
\frac{1}{K}|G/M_k|\leq |H/N_{\alpha(k)}|\leq K.|G/M_k|
\end{equation}
Set $r_k=\log_2|G/M_k|$ and $s_k=\log_2|H/N_k|$, so that $|r_k-s_{\alpha(k)}|\leq \log_2K$. Observe that, as $(M_k)_{k>0}$ is a strict filtration, we have $[M_k:M_{k+1}]\geq 2$ hence, for $\ell\geq k$, we have $|G/M_\ell|\geq 2^{\ell-k}|G/M_k|$ and $r_\ell-r_k\geq\ell-k$. We now prove that $\alpha$ satisfies the assumption of Lemma \ref{permut}. Let $N$ be an integer such that $N>2\log_2K$, so that the intervals $[r_k-\log_2K,r_k+\log_2K]$ and $[r_\ell-\log_2K,r_\ell+\log_2K]$ are disjoint if $\ell\geq k+N$. So, if $\ell\geq k+N$, then $s_{\alpha(k)}<s_{\alpha(\ell)}$, hence $\alpha(k)<\alpha(\ell)$ as $(N_k)_{k>0}$ is a strict filtration. The result now follows from Lemma \ref{permut}.
\hfill$\square$

\begin{Rem}
The above is not true in the more general setting of Lemma \ref{main}, i.e. simply considering coarse disjoint unions of graphs. Indeed, letting $X_k=Y_k$ be the $k$-cycle graph, one can construct a coarse equivalence between $X=\coprod_{k>0} X_k$ and $Y=\coprod_{k>0} Y_k$ as follows. First, group the connected components as $$\{X_1\}, \{X_2, X_3\}, \{X_4,X_5,X_6\}, \{X_7, X_8, X_9, X_{10}\},...$$ Then consider a cyclic permutation of the indices of each subset, for example, $4\mapsto 5, 5 \mapsto 6, 6 \mapsto 4$. Use this to define a map between the connected components of $X$ and $Y$, i.e., $X_4\mapsto Y_5, X_5 \mapsto Y_6, X_6 \mapsto Y_4$, such that on each connected component this map is a coarse equivalence (we can choose this so that the control functions are uniform across all $k$). The resulting map is a coarse equivalence between $X$ and $Y$, but the bijection induced on the indices does not have bounded displacement for any cofinite subset of $\mathbb{N}$. Note that $X$ and $Y$ are in fact both equal to the full box space of $\mathbb{Z}$. 
\end{Rem}

The above results give a necessary condition for box spaces to be coarsely equivalent.

\begin{Cor}\label{necessary} Let $(M_k)_{k>0}, (N_k)_{k>0}$ be strict filtrations of the residually finite groups $G,H$ respectively. If the box spaces $\square_{(M_k)}G$ and $\square_{(N_k)}H$ are coarsely equivalent, then there exists an almost permutation $\alpha$ of $\N$, with bounded displacement, such that the sequences $(\frac{|G/M_k|}{|H/N_{\alpha(k)}|})_{k\gg 0}$ and $(\frac{|H/N_{\alpha(k)}|}{|G/M_k|})_{k\gg 0}$ are bounded.
\end{Cor}

{\bf Proof:} This follows from the double inequality (\ref{1symm}), combined with the statement of Proposition \ref{bdddispl}.
\hfill$\square$

\begin{Lem}\label{uncountable} For $s$ a real number, $s\geq 1$, set $N_k(s)=2^{[ks]}$, where $[.]$ denotes the floor function. For $s,t\geq 1$, there exists an almost permutation $\alpha$ of $\N$ with bounded displacement such that the sequences $(\frac{N_k(s)}{N_{\alpha(k)}(t)})_{k\gg 0}$ and $(\frac{N_{\alpha(k)}(t)}{N_k(s)})_{k\gg 0}$ are bounded, if and only if $s=t$.
\end{Lem}

{\bf Proof:} Assume that the almost permutation $\alpha$ exists. Then there exists $K>0$ such that $|[ks]-[\alpha(k)t]|\leq K$ for $k\gg 0$. Observe that $\lim_{k\rightarrow\infty}\frac{[ks]}{k}=s$, so $\lim_{k\rightarrow\infty}\frac{[\alpha(k)t]}{k}=s$. Let $N\in\N$ be such that $k-N\leq\alpha(k)\leq k+N$ for every $k\in\N$. Then:
$$\frac{[\alpha(k)t]}{k+N}\leq\frac{[\alpha(k)t]}{\alpha(k)}\leq\frac{[\alpha(k)t]}{k-N}.$$
For $k\rightarrow\infty$, we deduce $s=\lim_{k\rightarrow\infty}\frac{[\alpha(k)t]}{\alpha(k)}$, but the latter limit is clearly equal to $t$.
\hfill$\square$

\medskip
For $\Z$, the converse of Corollary \ref{necessary} holds.

\begin{Cor} Let $(M_k)_{k>0}, (N_k)_{k>0}$ be two increasing sequences of positive integers, with $M_k$ dividing $M_{k+1}$ and $N_k$ dividing $N_{k+1}$. 
\begin{enumerate}
\item[1)] The box spaces $\square_{(M_k\Z)}\Z$ and $\square_{(N_k\Z)}\Z$ are coarsely equivalent if and only if there exists an almost permutation $\alpha$ of $\N$, with bounded displacement, such that the sequences $(\frac{M_k}{N_{\alpha(k)}})_{k\gg 0}$ and $(\frac{N_{\alpha(k)}}{M_k})_{k\gg 0}$ are bounded.
\item[2)] For $1\leq s$, the box spaces $\square_{(N_k(s)\Z)}\Z$ are pairwise not coarsely equivalent.
\end{enumerate}
\end{Cor}

{\bf Proof:} The first statement follows from the fact that cycles $C_M$ and $C_N$ are quasi-isometric by an $(A,A)$-quasi-isometry, with $A=1+[\frac{N}{M}]$. The second statement follows from Lemma \ref{uncountable}.
\hfill$\square$

\medskip
Corollary \ref{necessary} is enough to distinguish some expanders up to coarse equivalence. For $m\geq 2$, set $\Gamma_m=SL_m(\Z)$. For $N\geq 2$, let $\Gamma_m(N)=:\ker[SL_m(\Z)\rightarrow SL_m(\Z/N\Z)]$ be the congruence subgroup. Theorem \ref{SLm(pk)} is contained in the next result. Note that for $m\neq n$, the first statement of Proposition \ref{SLm} is contained in Theorem \ref{Hilbert} (proved in the next section), which however rests on much deeper results.

\begin{Prop}\label{SLm} Let $m,n\geq 2$.
\begin{enumerate}
\item[1)] Let $p,q$ be primes. The box spaces $\square_{\Gamma_m(p^k)}\Gamma_m$ and $\square_{\Gamma_n(q^k)}\Gamma_n$ are coarsely equivalent if and only if they are equal, i.e. $m=n$ and $p=q$.
\item[2)] For $s\geq 1$, the box spaces $\square_{(\Gamma_m(N_k(s)))}\Gamma_m$ are pairwise not coarsely equivalent. For $m\geq 3$, they have the geometric property (T) of Willett-Yu \cite{WiYu}.
\end{enumerate}
\end{Prop}

{\bf Proof:} 
\begin{enumerate}
\item[1)]Assume that $\square_{\Gamma_m(p^k)}\Gamma_n$ and $\square_{\Gamma_n(q^k)}\Gamma_n$ are coarsely equivalent. By Corollary \ref{necessary} there exists an almost permutation $\alpha$ of $\N$ with bounded displacement, such that the sequences $(\frac{\Gamma_m/\Gamma_m(p^k)}{\Gamma_n/\Gamma_n(q^{\alpha(k)})})_{k\gg 0}$ and $(\frac{\Gamma_n/\Gamma_n(q^{\alpha(k)})}    {\Gamma_m/\Gamma_m(p^k)})_{k\gg 0}$ are bounded. Now:
$$|\Gamma_m/\Gamma_m(p^k)|=|SL_m(\Z/p^k\Z)|=p^{k(m^2-1)}(1-\frac{1}{p^m})(1-\frac{1}{p^{m-1}})...(1-\frac{1}{p^2}),$$
so the assumption implies that the sequences $(\frac{p^{k(m^2-1)}}{q^{\alpha(k)(n^2-1)}})_{k\gg 0}$ and $(\frac{q^{\alpha(k)(n^2-1)}}{p^{k(m^2-1)}})_{k\gg 0}$ are bounded. Let $N>0$ be such that $|\alpha(k)-k|\leq N$ for every $k$. From the inequalities:
$$\frac{p^{k(m^2-1)}}{q^{(k+N)(n^2-1)}}\leq \frac{p^{k(m^2-1)}}{q^{\alpha(k)(n^2-1)}}\leq \frac{p^{k(m^2-1)}}{q^{(k-N)(n^2-1)}},$$
we get that the sequences $(\frac{p^{k(m^2-1)}}{q^{k(n^2-1)}})_{k\gg 0}$ and $\frac{q^{k(n^2-1)}}{p^{k(m^2-1)}})_{k\gg 0}$ are bounded. This implies $p^{m^2-1}=q^{n^2-1}$, i.e. $p=q$ and $m=n$.
\item[2)] Proceeding as in (1) and appealing again to Corollary \ref{necessary}, we see that if $\square_{(\Gamma_m(N_k(s)))}\Gamma_m$ and $\square_{(\Gamma_m(N_k(t)))}\Gamma_m$ are coarsely equivalent, then there exists an almost permutation $\alpha$ of $\N$, with bounded displacement, such that the sequences $((\frac{N_k(s)}{N_k(t)})^{m^2-1})_{k\gg 0}$ and $((\frac{N_k(t)}{N_k(s)})^{m^2-1})_{k\gg 0}$ are bounded. By Lemma \ref{uncountable}, this implies $s=t$. Since $\Gamma_m$ has property (T) for $m\geq 3$, by Theorem 1.1(4) in \cite{WiYu}, any box space of $\Gamma_m$ has geometric property (T).
\end{enumerate}
\hfill$\square$

\section{Coarse equivalence and quasi-isometry}

Since our goal is to distinguish expanders with bounded girth up to coarse equivalence, we briefly recall two results that are certainly well-known to experts. Although having unbounded girth can be used (see \cite{MeNa, Hume}) to distinguish expanders, it is clear that it is not a coarse invariant, as one may add a few edges to each component to create short loops. For box spaces, we have:

\begin{Lem}\label{largegirth} TFAE:
\begin{itemize}
\item[i)] $\square_{(N_k),S} G$ has unbounded girth;
\item[ii)] The Cayley graph $\mathcal{G}(G,S)$ is a tree.
\end{itemize}
\end{Lem}

Note that $G$ admits a Cayley graph which is a tree, if and only if $G$ is isomorphic to $\mathbb{F}_k\star(\star_1^m \mathbb{Z}/2\mathbb{Z})$ (see \cite{FTN}).
\medskip

{\bf Proof of Lemma \ref{largegirth}:} 

$(i)\Rightarrow(ii)$ If $\mathcal{G}(G,S)$ is not a tree, it contains a loop of length $n$, so $g(G/N)\leq n$ for every $N$.

$(ii)\Rightarrow (i)$ If $\mathcal{G}(G,S)$ is a tree, for $N\triangleleft G$ of finite index, the girth of $G/N$ is the distance from $1$ to $N-\{1\}$ in $\mathcal{G}(G,S)$. By residual finiteness, this distance tends to infinity.
\hfill$\square$

\medskip
For expanders, we observe:

\begin{Lem} Being an expander is a coarse invariant of box spaces.
\end{Lem}

{\bf Sketch of proof:} It is standard that $\square_{(M_k)} G$ is an expander if and only if the non-connected graph $\coprod_k G/M_k$ satisfies a strong isoperimetric inequality. By Lemma \ref{main}, a coarse equivalence $f:\square_{(M_k)} G\rightarrow\square_{(N_k)} H$ is a quasi-isometry with uniform constants in restriction to any component $G/M_k$, so one can copy the proof that satisfying a strong isoperimetric inequality is invariant under quasi-isometry (as given e.g. in Theorem 7.34 of \cite{Soardi}).
\hfill$\square$

\medskip
A {\bf marked group} on $m$ generators is a group given as a quotient of the free group $\F_m$. We denote by $\mathcal{G}_m$ the space of marked groups, i.e. the set of quotients of $\F_m$, equipped with the Chabauty topology, that turns it into a compact, totally disconnected space. Recall that we denote by $B_G(N)$ the ball of radius $N$ centered at the identity of the marked group $G$.

\begin{Prop}\label{marked} Let $(G_k)_{k>0}$ be a sequence in $\mathcal{G}_m$ converging to $G\in\mathcal{G}_m$; similarly, let $(H_k)_{k>0}$ be a sequence in $\mathcal{G}_n$ converging to $H\in\mathcal{G}_n$. Assume that, for some $A>0$, there exists for every $k>0$ an $(A,A)$-quasi-isometry $f_k:G_k\rightarrow H_k$. Then there exists an $(A,A)$-quasi-isometry $f:G\rightarrow H$.
\end{Prop}

{\bf Proof:} Without loss of generality, we may assume that $f_k$ maps the unit of $G_k$ to the unit of $H_k$. Fix some radius $N>0$. For $k\gg 0$, the balls $B_G(N)$ and $B_{G_k}(N)$ (resp. $B_H(N)$ and $B_{H_k}(N)$ are isomorphic as labeled graphs, let $\phi_{k,N}: B_G(N)\rightarrow B_{G_k}(N)$ (resp. $\psi_{k,N}:B_H(N)\rightarrow B_{H_k}(N)$) be a labeled isomorphism. Now $f_k$ maps $B_{G_k}(N)$ into $B_{H_k}(A(N+1))$ so we may define maps $\tilde{f}_{k,N}:B_G(N)\rightarrow B_H(A(N+1))$ as:
$$\tilde{f}_{k,N}: \psi^{-1}_{k,A(N+1)}\circ f_k\circ \phi_k.$$
Set $f(e)=e$. Observe that there are finitely many maps $B_G(N)\rightarrow B_H(A(N+1))$, so we may set up a diagonal process: first find an infinite set $I_1$ of indices $k$ on which $\tilde{f}_{1,k}$ is constant, define $\tilde{f}_1=\tilde{f}_{1,k}$ for $k\in I_1$; then find an infinite subset $I_2\subset I_1$ on which $\tilde{f}_{2,k}$ is constant, define $\tilde{f}_2=\tilde{f}_{2,k}$ for $k\in I_2$, so that $\tilde{f}_2$ restricts to $\tilde{f}_1$ on $B_G(1)$; iterating this construction, we obtain a well-defined $f:G\rightarrow H$ whose restriction to $B_G(N)$ is $\tilde{f}_N$.
\begin{itemize}
\item $f$ is an $(A,A)$-quasi-isometric embedding: for $x,y\in G$, take $N$ large enough that $B_G(N)$ contains $x^{-1}y$ and $B_H(A(N+1))$ contains $f(x)^{-1}f(y)$; then $d_H(f(x),f(y))=d_{H_k}(f_k(\phi_{k,N}(x)),f_k(\phi_{k,N}(y)))$ for $k\in I_N$, so 
$$\frac{1}{A}d_{G_k}(\phi_{k,N}(x),\phi_{k,N}(y))-A\leq d_H(f(x),f(y))\leq A.d_{G_k}(\phi_{k,N}(x),\phi_{k,N}(y)) +A,$$
but $d_{G_k}(\phi_{k,N}(x),\phi_{k,N}(y))=d_G(x,y)$.
\item Let $C>0$ be such that $H_k$ is the $C$-neighbourhood of $f_k(G_k)$ for every $k$. We claim that $H$ is the $C$-neighbourhood of $f(G)$. So, fix $y\in H$. Take $N\geq A(A+C+d_H(1,y))$. Choose $k\in I_N$. Since $H_k$ is the $C$-neighbourhood of $f_k(G_k)$, we find $x\in G$ such that $d_{H_k}(\psi_{k,N}(y),f_k(\phi_{k,N}(x)))\leq C$. By the triangle inequality, and using $f_k(e)=e$:
$$d_H(1,y)+C\geq d_{H_k}(1,f_k(\phi_{k,N}(x)))\geq \frac{1}{A}d_{G_k}(1,\phi_{k,N}(x))-A,$$
hence $d_{G_k}(1,\phi_{k,N}(x))\leq N$. This shows that $x\in B_G(N)$ and therefore $\tilde{f}_{N,k}(x)$ is independent of the choice of $k\in I_N$. Finally $d_H(f(x),y)=d_{H_k}(f_k(\phi_{k,N}(x)),\psi_{k,N}(y))\leq C$.
\end{itemize}\hfill$\square$

The following result was suggested to us by R. Tessera: 

\begin{Thm}\label{coarse=>qi}\footnote{During the writing process of this paper we were informed that K. Das obtained the stronger result that, under the same assumptions, $G$ and $H$ are uniformly measure equivalent.} Let $G,H$ be residually finite, finitely generated groups. If $\square_{(M_k)} G$ is coarsely equivalent to $\square_{(N_k)} H$, then $G$ is quasi-isometric to $H$. 
\end{Thm}

{\bf Proof:} Since the sequence $(G/M_k)_{k>0}$ converges to $G$ and the sequence $(H/N_k)_{k>0}$ converges to $H$ in the space of marked groups, the result follows by combining Lemma \ref{main} and Proposition \ref{marked}.
\hfill$\square$

\medskip
{\bf Proof of Theorem \ref{Hilbert}:} 

(1) and (2). Since $SL_2(\Z[\sqrt{p}])$ has property $(\tau)$ (see section 4.3.3 in \cite{Lub}), $X_p$ is an expander. Similarly, since $SL_n(\Z)$ has Kazhdan's property (T), $Y_n$ is an expander. View $SL_2(\Z[\sqrt{p}])$ (resp. $SL_n(\Z)$) as a lattice in $SL_2(\R)\times SL_2(\R)$ (resp. $SL_n(\R)$). If $X_p$ is coarsely equivalent to $X_q$ (resp. $Y_m$ is coarsely equivalent to $Y_n$) then, by Theorem \ref{coarse=>qi}, $SL_2(\Z[\sqrt{p}])$ is quasi-isometric to $SL_2(\Z[\sqrt{q}])$ (resp. $SL_m(\Z)$ is quasi-isometric to $SL_n(\Z)$). By quasi-isometric rigidity of lattices (see \cite{FS} for $SL_2(\Z[\sqrt{p}])$, and \cite{Eskin} for $SL_n(\Z)$) this implies $p=q$ (resp. $m=n$).

(3) Since $SL_2(\Z[\sqrt{p}])$ is a-(T)-menable, $X_p$ admits a fibered coarse embedding into Hilbert space (see \cite{ChWaWa}, especially Theorem 1.1); since $SL_n(\Z)$ has Kazhdan's property (T), $Y_n$ has geometric property (T) (see \cite{WiYu}, especially Theorem 1.1(4)). So, if there were a coarse embedding of $Y_n$ in $X_p$, then $Y_n$ admits a fibered coarse embedding into Hilbert space, and this contradicts Theorem 8.1(2) in \cite{WiYu}.\footnote{According to Example 2.6 in \cite{ChWaWa}, a sequence of graphs with unbounded girth admits a fibered coarse embedding into Hilbert space, so  by the same argument there is no coarse embedding of $Y_n$ into such a sequence of graphs - in particular into the expanders from \cite{Hume}.}
\hfill$\square$

\medskip
If $(X_n)_{n>0}$ is a sequence of bounded metric spaces, we define the {\bf coarse disjoint union} as $\coprod_{n=1}^\infty X_n$ endowed (as we did for graphs in the introduction) with a metric $d$ inducing the original one on each $X_n$ and such that $d(X_m,X_n)\geq \max\{diam(X_m),diam(X_n)\}$ for $m\neq n$.

\begin{Cor}\label{buildings} Let $X,Y$ be either irreducible Riemannian symmetric spaces of rank $\geq 2$, or irreducible Euclidean buildings of rank $\geq 2$. Let $\Gamma,\Lambda$ be co-compact lattices on $X, Y$ respectively. Fix a filtration $(M_k)_{k>0}$ in $\Gamma$ (resp. $(N_k)_{k>0}$ in $\Lambda$). Assume that the coarse disjoint unions $\coprod_{k>0} M_k\backslash X$ and $\coprod_{k>0} N_k\backslash Y$ are coarsely equivalent. Then $X$ is isometric to $Y$.
\end{Cor}

{\bf Proof:}
\begin{enumerate}
\item By the Milnor-\v{S}varc lemma, any orbit map $\Gamma\rightarrow X$ induces a quasi-isometry between $\Gamma$ and $X$; this map being $M_k$-equivariant for every $k>0$, it descends to a quasi-isometry between the box space $\square_{(M_k)} \Gamma$ and the coarse union $\coprod_{k>0} M_k\backslash X$; and similarly for $\square_{(N_k)}\Lambda$ and $\coprod_{k>0} N_k\backslash Y$.
\item Since $\square_{(M_k)}\Gamma$ and $\square_{(N_k)}\Lambda$ are coarsely equivalent, by Theorem \ref{coarse=>qi}, $\Gamma$ and $\Lambda$ are quasi-isometric, hence $X$ and $Y$ are quasi-isometric. 
\item By a celebrated result of B. Kleiner and B. Leeb (Theorem 1.1.3 in \cite{KlLe}), $X$ and $Y$ are isometric. 
\end{enumerate}\hfill$\square$

\medskip
This allows us to answer a question asked by A. Lubotzky (personal communication): in Theorem 1.2 of \cite{LSV}, $(d-1)$-dimensional Ramanujan complexes are defined as quotients of the Euclidean building of $PGL_d(F)$, where $F$ is a local field of positive characteristic, by congruence subgroups of a suitably defined arithmetic subgroup. Consider now two families of Ramanujan complexes, one from $PGL_d(F)$, one from $PGL_{d'}(F')$, with $d,d'\geq 3$. Corollary \ref{buildings} implies that, if the coarse disjoint unions of these families are coarsely equivalent, then $d=d'$ and $F=F'$. Of course the restriction to $d\geq 3$ leaves open the very interesting question whether the families of $(p+1)$-regular Ramanujan graphs constructed in \cite{LPS} are coarsely equivalent or not, for distinct primes $p$.

\section{Large diameter}

\subsection{Around property $D_\alpha$}

In this section we consider property $D_\alpha$ (with $0<\alpha\leq 1$). Our first result is that $D_\alpha$ is coarsely invariant.

\begin{Lem} Property $D_\alpha$ is a coarse invariant of box spaces.
\end{Lem}

{\bf Proof:} Let $f:\square_{(M_k)} G\rightarrow\square_{(N_k)} H$ be a coarse equivalence, with $\square_{(M_k)} G$ satisfying $D_\alpha$. Let $H/N$ be a component of $\square_{(N_k)} H$. By Lemma \ref{main}, we may assume that $H/N$ contains the image by $f$ of a unique component $G/M$ of $\square_{(M_k)} G$. Then:
$$diam(H/N)\geq diam(f(G/M))\geq \frac{1}{A}diam(G/M)-A\geq \frac{K}{A}.|G/M|^\alpha-A$$ 
By equation (\ref{volume}):
$$diam(H/N)\geq\frac{K}{A.|B_H(C)|^\alpha} |H/N|^\alpha - A$$
hence $\square_{(N_k)} H$ has property $D_\alpha$.
\hfill$\square$


\medskip
As mentioned in the introduction, it was proved in \cite{BrTo} that if $G$ admits a box space with property $D_\alpha$, then $G$ virtually maps onto $\Z$. Theorem \ref{Breuillard} then follows from the next result. We thank E. Breuillard for suggesting working with $G_0/[G_0,G_0]$.

\begin{Prop}\label{virtuallymapping} Let $G$ be a (finitely generated, residually finite) group virtually mapping onto $\Z$. Let $G_0$ be a finite index, normal subgroup of $G$ with infinite abelianization. Let $d\geq 1$ be the rank of the free abelian part of $G_0/[G_0,G_0]$ (i.e. $d=\dim_\Q (G_0/[G_0,G_0])\otimes_\Z \Q$). Then, for any $\alpha<\frac{1}{d}$, the group $G$ admits a box space with property $D_\alpha$.
\end{Prop}

{\bf Proof:} Let $(M_n)_{n>0}$ be a filtration of $G$; let $\epsilon>0$ be defined by $\alpha=\frac{1}{d+\epsilon}$.

For $k>0$, we define $L_k=G_0^k[G_0,G_0]$, the verbal subgroup of $G_0$ generated by all elements of the form $a^k[b,c]\;(a,b,c\in G_0)$: it is a finite index, characteristic subgroup of $G_0$, so it is a finite index, normal subgroup of $G$. 

Let $(k_n)_{n>0}$ be a sequence of positive integers, with $k_n$ dividing $k_{n+1}$ and $k_n^\epsilon \geq |G/M_n|$ for every $n>0$. Set $N_n=:L_{k_n}\cap M_n$. The sequence $(N_n)_{n>0}$ is a filtration of $G$, and we claim that the box space $\square_{(N_n)} G$ has property $D_\alpha$.

To prove this, we exploit the short exact sequence
$$1\rightarrow L_{k_n}/N_n\rightarrow G/N_n\rightarrow G/L_{k_n}\rightarrow 1.$$
Since $G_0/L_{k_n}$ maps onto the cyclic group of order $k_n$, we have
\begin{equation}\label{(3)}
diam(G/N_n)\geq diam(G/L_{k_n})\geq\frac{1}{[G:G_0]}diam(G_0/L_{k_n})\geq\frac{k_n}{2[G:G_0]}.
\end{equation}
On the other hand, since $G/L_{k_n}$ has a subgroup of bounded index isomorphic to $(\Z/k_n\Z)^d$, there exists a constant $A>0$ such that
$$|G/L_{k_n}|\leq A.k_n^d$$
Moreover:
$$|L_{k_n}/N_n|=|L_{k_n}/(L_{k_n}\cap M_n)|=|L_{k_n}M_n/M_n|\leq|G/M_n|\leq k_n^\epsilon.$$
So
\begin{equation}\label{(4)}
|G/N_n|=|G/L_{k_n}|.|L_{k_n}/N_n|\leq A.k_n^{d+\epsilon}.
\end{equation}
Comparing (\ref{(3)}) and (\ref{(4)}) gives property $D_\alpha$ with $\alpha=\frac{1}{d+\epsilon}$.
\hfill$\square$

\medskip
We offer a different construction for groups mapping onto $\Z$. Our excuse for presenting it is that it gives a better range of $\alpha$'s than Proposition \ref{virtuallymapping}.

\begin{Prop}\label{MappingontoZ} Assume that $G$ maps onto $\Z$. Fix $\alpha\in]0,1[$. There exists a box space of $G$ with property $D_\alpha$.
\end{Prop}

{\bf Proof:} Write $G$ as a semi-direct product $G=H\rtimes\Z$, with $H$ not necessarily finitely generated; denote by $t\in G$ a generator of the factor $\Z$. Choose a finite generating set of $G$ of the form $S=\{t\}\cup T$, with $T\subset H$. Let $(M_n)_{n>0}$ be a filtration of $G$. Define inductively a sequence $(k_n)_{n>0}$ of positive integers satisfying:
\begin{itemize}
\item $k_{n-1}$ divides $k_n$;
\item for every $x\in H$, we have $[t^{k_n},x]\in H\cap M_n$ (i.e. $k_n$ is a multiple of the order of $Ad(t)$ on $H/(H\cap M_n)$).
\item $k_n\geq 2^{\frac{1}{1-\alpha}}|H/(H\cap M_n)|^{\frac{\alpha}{1-\alpha}}.$
\end{itemize}
Set then $N_n=:<H\cap M_n,t^{k_n}>=\{at^{\ell k_n}:a\in H\cap M_n,\ell\in\Z\}$, a normal subgroup of $G$. The family $(N_n)_{n>0}$ is a filtration of $G$ and, for the corresponding box space $\square_{(N_n)} G$, we have $|G/N_n|=|H/(H\cap M_n)|.k_n$ and
$$diam(G/N_n)\geq \frac{k_n}{2}\geq |G/N_n|^\alpha$$
by choice of $k_n$. This is property $D_\alpha$.
\hfill$\square$

\medskip
The following result came out of a conversation with V. Alekseev. It shows that any box space can be coarsely embedded in a box space with property $D_\alpha$.

\begin{Prop} Assume that $\square_{(M_k)} G$ has property $D_\alpha$. Let $\epsilon$ be such that $0<\epsilon<\alpha$. Then, for any finitely generated residually finite group $H$, and any filtration $(N_k)_{k>0}$ of $H$, there exists a sub-filtration $(M_{(\ell_k)})$ of $(M_k)$ such that $\square_{(M_{\ell_k}\times N_k)}(G\times H)$ has property $D_{\alpha-\epsilon}$.
\end{Prop}

{\bf Proof:} Select $\ell_k$ such that $|G/M_{\ell_k}|\geq |H/N_k|^{\frac{\alpha}{\epsilon}-1}$. Then $|G/M_{\ell_k}|^\alpha\geq (|G/M_{\ell_k}|.|H/N_k|)^{\alpha-\epsilon}$, hence by property $D_\alpha$ for $\square_{(M_{\ell_k})} G$:
$$diam((G\times H)/(M_{\ell_k}\times N_k))\geq diam(G/M_{\ell_k})\geq C.|G/M_{\ell_k}|^\alpha\geq C.(|G/M_{\ell_k}|.|H/N_k|)^{\alpha-\epsilon}$$
which is the announced result.
\hfill$\square$

\medskip
We now prove Theorem \ref{virtuallycyclic} from the introduction.

\medskip
{\bf Proof of Theorem \ref{virtuallycyclic}:} $(i)\Rightarrow (ii)$ Suppose that $diam(G/M)\geq K.|G/M|$ for every component $G/M$ of $\square G$. We will show that $G$ has linear growth, hence is virtually cyclic (see Theorem 3.1 in \cite{Mann}; note that this result does NOT appeal to Gromov's polynomial growth theorem). So fix $R>0$. Let $x,y\in G/M$ be such that $d(x,y)=diam(G/M)$, and let $k$ be the integer part of $\frac{diam(G/M)}{2R+1}$. On some geodesic from $x$ to $y$, place $k$ points $x_1,...,x_k$ such that $d(x_i,x_{i+1})=2R+1$. Then the balls $B_{G/M}(x_i,R)$ are pairwise disjoint with same cardinality $B_{G/M}(R)$, so we get:
$$k.B_{G/M}(R)\leq |G/M|\leq\frac{1}{K}.diam(G/M)\leq\frac{k}{K}.(2R+2),$$
i.e. $B_{G/M}(R)\leq \frac{1}{K}.(2R+2)$. Now use residual finiteness to find an $M$ such that the ball of $G$ of radius $R$ centered at the identity, injects into $G/M$. For such an $M$ we may replace $B_{G/M}(R)$ by $B_G(R)$, establishing that $G$ has linear growth.

\medskip
$(ii)\Rightarrow(i)$ If $G$ is virtually cyclic, let $H\triangleleft G$ be an infinite, cyclic subgroup of finite index $d$. The inclusion $H\subset G$ is a quasi-isometry which moreover is right-$H$-equivariant. So, for $M$ running along finite index normal subgroups of $G$, it descends to a family of quasi-isometries (with uniform constants!) between $H/H\cap M\subset G/M$ and $G/M$. So the diameters of $H/H\cap M$ and $G/M$ are the same, up to a multiplicative constant not depending on $M$. Now $H/H\cap M$ is a cycle, so its diameter is one half of $|H/H\cap M|$. Finally we have 
$$|H/H\cap M|\leq |G/M|\leq d.|H/H\cap M|,$$
concluding the proof.
\hfill$\square$

\subsection{The lamplighter group}
We discuss more closely the lamplighter group $G=(\Z/2\Z)\wr\Z$ (a solvable,  non-polycyclic group). It was observed in Remark 5.3 of \cite{BrTo} that finite lamplighter groups $G_n=(\Z/2\Z)\wr(\Z/n\Z)$ are not almost flat, as $diam(G_n)=O(n)$. This shows that the coarse disjoint union of the $G_{2^k}$'s, viewed as a box space of $G$, does not have property $D_\alpha$ for any $\alpha$. However, by Proposition \ref{MappingontoZ} above, $G$ admits box spaces with property $D_\alpha$. We construct explicitly such a box space.



We preface the construction with a number of algebraic remarks. We denote by $F_{2^n}$ the finite field of order $2^n$, and by $p_i$ the $i$-th prime. Let $\alpha_i$ be a generator of the multiplicative group of $F_{2^{p_i}}$. Let $P_i(X)\in F_2[X]$ be the minimal polynomial of $\alpha_i$ over $F_2$ (so $\deg P_i=p_i$), and $(P_i)$ the principal ideal of $F_2[X]$ generated by $P_i$. Then
$$F_2[X]/(P_1)...(P_k)\simeq \prod_{i=1}^k F_2[X]/(P_i),$$
so at the level of multiplicative groups:
$$(F_2[X]/(P_1)...(P_k))^\times \simeq \prod_{i=1}^k (F_2[X]/(P_i))^\times.$$
Since $2^{p_1}-1,...,2^{p_k}-1$ are pairwise relatively prime, by the Chinese remainder theorem $\prod_{i=1}^k (F_2[X]/(P_i))^\times$ is cyclic of order $\ell_k=:\prod_{i=1}^k (2^{p_i}-1)$, and the class of $X$ in $(F_2[X]/(P_1)...(P_k))^\times$ has order $\ell_k$. Observe that $\prod_{i=1}^k (2^{p_i}-1)\geq C.\prod_{i=1}^k 2^{p_i}$ for some universal constant $C>0$: indeed the infinite product $\prod_{i=1}^\infty (1-2^{-p_i})$ converges to  a positive number, as $\sum_{i=1}^\infty 2^{-p_i}\leq \sum_{i=1}^\infty 2^{-i}=1$.

Set then $A=F_2[X,X^{-1}]$, the ring of Laurent polynomials with coefficients in $F_2$; we view $A$ as $F_2[X]$ localized with respect to $X$. We then identify $G$ with the following subgroup of $GL_2(A)$:
$$G=\{\left(\begin{array}{cc}X^n & P \\0 & 1\end{array}\right):P\in A, n\in\Z\}.$$
As a generating set of $G$, we take $\{\left(\begin{array}{cc}X^{\pm 1} & 0 \\0 & 1\end{array}\right), \left(\begin{array}{cc}1 & 1 \\0 & 1\end{array}\right)\}$. Let $I_k$ be the ideal of $A$ generated by the product $P_1...P_k$, so that $A/I_k\simeq F_2[X]/(P_1)...(P_k)$. Let $N_k$ be the intersection of $G$ with the congruence subgroup $\ker(GL_2(A)\rightarrow GL_2(A/I_k))$ associated with $I_k$. The family $(N_k)_{k>0}$ is a filtration of $G$.

\begin{Prop}\label{lamplighter} For $G=(\Z/2\Z)\wr\Z$, the above box space $\square_{(N_k)}$ has property $D_{1/2}$.
\end{Prop}

{\bf Proof:} We have $|G/N_k|=\ell_k.\prod_{i=1}^k 2^{p_i}$ and
$$diam(G/N_k)\geq diam(C_{\ell_k})=\frac{\ell_k}{2}\geq\frac{C^{1/2}}{2}|G/N_k|^{1/2}$$
which is property $D_{1/2}$.
\hfill$\square$


\subsection{A lattice in SOL}

We consider an explicit polycyclic group, namely the semi-direct product $\Gamma=\Z^2\rtimes\Z$, where $\Z$ acts by powers of $A=\left(\begin{array}{cc}1 & 1 \\1 & 0\end{array}\right)$; it is a uniform lattice in the 3-dimensional solvable Lie group SOL. We embed $\Gamma$ in a standard way in $GL_3(\Z)$. For $N>0$, we denote by $\Gamma(N)$ the {\it congruence subgroup}, i.e. the intersection of $\Gamma$ with the kernel of reduction modulo $N$:
$$\Gamma(N)=\Gamma\cap \ker(GL_3(\Z)\rightarrow GL_3(\Z/N\Z)).$$
If $(N_k)_{k>0}$ is an increasing sequence of integers, with $N_k$ dividing $N_{k+1}$, then $(\Gamma(N_k))_{k>0}$ is a filtration of $\Gamma$, and the corresponding box space $\square_{(\Gamma(N_k))}\Gamma$ is a {\it congruence box space}.

\begin{Prop}\label{SOL} 
\begin{enumerate}
\item[1)] The congruence box space $\square_{(\Gamma(5^k))}\Gamma$ has property $D_{1/3}$.
\item[2)] If a congruence box space $\square_{(\Gamma(N_k))}\Gamma$ has property $D_\alpha$, then $\alpha\leq\frac{1}{3}$.
\item[3)] There exists a congruence box space $\square_{(\Gamma(N_k))}\Gamma$ without property $D_\alpha$, for every $\alpha\in]0,1]$.
\end{enumerate}
\end{Prop}

So, for $\alpha$ close to 1, the box spaces of $\Gamma$ with property $D_\alpha$ constructed in Proposition \ref{MappingontoZ} are not coarsely equivalent to any congruence box space of $\Gamma$.

We denote by $\delta(N)$ the order of the matrix $A$ in $GL_2(\Z/N\Z)$. Note that $\Gamma/\Gamma(N)=(\Z/N\Z)^2\rtimes\Z/\delta(N)\Z$. It is classical that $A^k=\left(\begin{array}{cc}F_{k+1} & F_k \\F_k & F_{k-1}\end{array}\right)$, where $F_k$ is the $k$-th Fibonacci number (normalized by $F_0=0,F_1=1$); we take advantage of the fact that, because of this connection, the function $\delta(N)$ has been well-studied (see e.g. \cite{Ren}).

\begin{Lem}\label{Fibo} There exists a constant $C>0$ such that $C.\log(N)\leq\delta(N)\leq 6N$.
\end{Lem}

{\bf Proof of Lemma \ref{Fibo}:} The second inequality is result A.7 in \cite{Ren}. For the lower bound, we appeal to result A.8 in \cite{Ren}: let $L_t=(\frac{1+\sqrt{5}}{2})^t + (\frac{1-\sqrt{5}}{2})^t$ be the $t$-th Lucas number; if $t$ is the largest $t$ such that $L_t\leq N$, then $\delta(N)\geq 2t$. The lemma follows.
\hfill$\square$

\medskip
As a generating set of $G$ take e.g. $$S=\{((\pm1,0),Id); ((0,\pm 1),Id); ((0,0),A^{\pm 1})\}.$$

\begin{Lem}\label{diameter} There are constants $C_1,C_2>0$ such that, for large enough $N$:
$$C_1\delta(N)\leq diam(\Gamma/\Gamma(N))\leq C_2\delta(N)$$
\end{Lem}


\medskip
{\bf Proof:} For the lower bound, observe that $\Gamma/\Gamma(N)$ maps onto the cycle $C_{\delta(N)}$, so $\frac{\delta(N)}{2}=diam(C_{\delta(N)})\leq diam (\Gamma/\Gamma(N))$. For the upper bound: let $((m,n),k)\in\Gamma/\Gamma(N)=(\Z/N\Z)^2\rtimes\Z/\delta(N)\Z$ that realizes the maximum of the word length (with $0\leq m,n<N;\;0\leq k<\delta(N)$). Since the quotient map $\Gamma\rightarrow\Gamma/\Gamma(N)$ decreases distances, we have $diam(\Gamma/\Gamma(N))\leq |((m,n),k)|_S$ (where in the RHS, $((m,n),k)$ is now viewed as an element of $\Gamma$). So by the triangle inequality 
$$diam(\Gamma/\Gamma(N))\leq |((m,n),0)|_S + |((0,0),k)|_S\leq |((m,n),0)|_S+ \frac{\delta(N)}{2}.$$
Since $\Z^2$ is exponentially distorted in $\Gamma$, there exists a constant $C>0$ such that the intrinsic and extrinsic word lengths satisfy $|((a,b),0)|_S\leq C.\log(1+|(a,b)|_{\Z^2})$ for every $(a,b)\in\Z^2$ (see e.g. Lemma 3.4 in \cite{Osin}). So we get:
$$diam(\Gamma/\Gamma(N))\leq C.\log(1+|(m,n)|_{\Z^2})+\frac{\delta(N)}{2}\leq C.\log(1+2N)+\frac{\delta(N)}{2}.$$
Combining with Lemma \ref{Fibo}, the lemma follows. 
\hfill$\square$

\medskip
{\bf Proof of Proposition \ref{SOL}:} 
\begin{enumerate}
\item[1)] Let $\alpha(N)$ be the smallest index so that $N$ divides $F_{\alpha(N)}$ (e.g. $\alpha(2)=3, \alpha(3)=4, \alpha(4)=6$, etc...). We clearly have $\delta(N)\geq\alpha(N)$; moreover $\alpha(5^k)=5^k$, by equality (35)  in \cite{Hal}. By Lemma \ref{diameter}, there is a constant $C>0$ such that:
\begin{equation}\label{5^k-1}
diam(\Gamma/\Gamma(5^k))\geq C.\delta(5^k)\geq C.5^k.
\end{equation}
On the other hand, by Lemma \ref{Fibo}:
\begin{equation}\label{5^k-2}
|\Gamma/\Gamma(5^k)|=5^{2k}.\delta(5^k)\leq 6.5^{3k}
\end{equation}
Comparing formulae (\ref{5^k-1}) and (\ref{5^k-2}) gives property $D_{1/3}$.

\item[2)] Assume that the congruence box space $\square_{(\Gamma(N_k))}\Gamma$ satisfies property $D_\alpha$. 
By Lemma \ref{diameter}, there exists $C>0$ such that:
\begin{equation}
diam(\Gamma/\Gamma(N_k))\leq C\delta(N_k).
\end{equation}
Plugging this into property $D_\alpha$:
$$\delta(N_k)\geq K.N_k^{2\alpha}\delta(N_k)^\alpha \Longleftrightarrow \delta(N_k)^{1-\alpha}\geq K.N_k^{2\alpha}.$$
By Lemma \ref{Fibo}:
$$(6N_k)^{1-\alpha}\geq K.N_k^{2\alpha}.$$
Since $N_k$ can be arbitrarily large, this forces $\alpha\leq\frac{1}{3}$.

\item[3)] We seek an increasing sequence $(N_k)_{k>0}$, with $N_k$ dividing $N_{k+1}$, such that for every $\alpha\in]0,1[$:
$$diam(\Gamma/\Gamma(N_k))=o(|\Gamma/\Gamma(N_k)|^\alpha)=o(N_k^{2\alpha}\delta(N_k)^\alpha).$$
In view of Lemma \ref{diameter}, it is enough to have:
\begin{equation}\label{nonDalpha}
\delta(N_k)^{1-\alpha}=o(N_k^{2\alpha}).
\end{equation}
We will use two properties of the Fibonacci numbers $F_n$:
\begin{itemize}
\item $\delta(F_n)=4n$ for odd $n$ (formula A.9 in \cite{Ren});
\item $F_m,F_n$ are relatively prime when $m,n$ are relatively prime (classical).
\end{itemize}
Let $q_i$ denote the $i$-th {\it odd} prime. Set $N_k=\prod_{i=1}^k F_{q_i}$. Since the $F_{q_i}$'s are pairwise relatively prime, we have, by the Chinese Remainder Theorem:
$$\delta(N_k)=\prod_{i=1}^k \delta(F_{q_i})=4^k.\prod_{i=1}^k q_i.$$
The function $\prod_{i=1}^k q_i$ is basically the {\it primorial} function $p_k\sharp$, i.e. the product of the first $k$ primes. It is known (see \cite{We1}) that $\lim_{k\rightarrow\infty}(p_k\sharp)^{1/p_k}=e$; using $q_k \sim k\log(k)$, this implies that there exists a constant $C>0$ with $\delta(N_k)\leq C^{k\log(k)}$.

On the other hand, let $K>1$ be such that $F_n\geq K^n$ for $n\geq 3$. Then 
$$N_k^{2\alpha}\geq K^{2\alpha\sum_{i=1}^k q_i}.$$
But $\sum_{i=1}^k q_i \sim k^2\log(k)$ (see \cite{We2}), from which (\ref{nonDalpha}) becomes clear.
\end{enumerate}
\hfill$\square$


\section{Full box spaces}

\subsection{Virtually cyclic groups}

We preface the proof of Theorem \ref{fullvirtuallycyclic} with a lemma:

\begin{Lem}\label{cyclic} Let $p:G\rightarrow Q$ with a surjective homomorphism with finite kernel $F$. Let $t\in Q$ be an element of infinite order. For $n\in\N$, the number $L_n$ of subgroups $H$ of $G$ such that $p(H)=<t^n>$, is bounded independently of $n$.
\end{Lem}

{\bf Proof:} If $p(H)=<t^n>$, let $h_0\in H$ be such that $p(h_0)=t^n$. Then any element $h\in H$ is uniquely of the form $h=h_0^k.f$, where $k\in \Z$ and $f\in F\cap H$. The subgroup $H$ is therefore determined by $h_0$ and the subgroup $F\cap H$; so a very rough upper bound for $L_n$ is $|F|\cdot 2^{|F|}$.
\hfill$\square$

\medskip
{\bf Proof of Theorem \ref{fullvirtuallycyclic}:} Let $G$ be a virtually cyclic group. It is a classical fact (see e.g. Lemma 11.4 in \cite{Hempel}) that $G$ maps with finite kernel either onto $\Z$ or onto the infinite dihedral group $D_\infty$. So we distinguish two cases.

\underline{1st case:} There is a surjective homomorphism $p:G\rightarrow\Z$ with finite kernel $F$. For $n\geq 1$, let $K_n$ be the number of normal subgroups $M\triangleleft G$ such that $p(M)=n\Z$.
We construct a quasi-isometry $f:\square G\rightarrow \square\Z$. If $M\triangleleft G$ is such that $p(M)=n\Z$, we first retract the Cayley graph of $G/M$ to the cycle $C_n$ corresponding to an element of order $n$ in $G/M$. In this way we get a quasi-isometry $\square G\rightarrow \coprod_{n\geq 1} K_n\cdot C_n$, where $K_n\cdot C_n$ is the disjoint union of $K_n$ cycles $C_n$. Second we construct a map $\coprod_{n\geq 1} C_n\rightarrow\coprod_{n\geq 1} K_n\cdot C_n$, by retracting $C_1,...,C_{K_1}$ onto $K_1\cdot C_1$, then retracting $C_{K_1+1},...,C_{K_1+K_2}$ onto $K_2\cdot C_2$, etc... This is a $(K,K)$-quasi-isometry, where $K=\max_n K_n$ (which is finite by Lemma \ref{cyclic}).

\underline{2nd case:} There is a surjective homomorphism $p:G\rightarrow D_\infty$ with finite kernel $F$. Recall that $D_\infty$ has a presentation $D_\infty=\langle s,t |s^2=t^2=1\rangle$. The group $D_\infty$ has 3 subgroups of index $2$ and, for $n\geq 2$, a unique subgroup $N_n$ of index $2n$; the latter is the kernel of the map $D_\infty\rightarrow D_n$, where $D_n$ is the dihedral group of order $2n$; so $N_n=<(st)^n>$. Let $K_n$ be the number of normal subgroups $M\triangleleft G$ with $p(M)=N_n$; proceeding as in the first step, we retract first $G/M$ onto $C_n$, then we built a map $\square\Z\rightarrow \coprod_{n\geq 1}K_n\cdot C_n$, using Lemma \ref{cyclic} to prove it is a quasi-isometry.
\hfill$\square$

\medskip
We will see in Corollary \ref{noncoarse} below, that Theorem \ref{fullvirtuallycyclic} cannot be extended to virtually-$\Z^2$ groups.

\subsection{Normal subgroup growth}

Denote by $a_n^\triangleleft(G)$ (resp. $s_n^\triangleleft(G)$) the number of normal subgroups of $G$ with index $n$ (resp. at most $n$). 

\begin{Prop}\label{subgroupgrowth} If $\square G$ and $\square H$ are coarsely equivalent, there exists constants $A,B>0$ such that, for $n\gg 0$:
$$a_n^\triangleleft(G)\leq \sum_{An\leq k\leq Bn} a_k^\triangleleft(H)=s_{Bn}^\triangleleft(H)-s_{An}^\triangleleft(H).$$
\hfill$\square$
\end{Prop}

{\bf Proof:} Let $f:\square G\rightarrow \square H$ be a coarse equivalence, with control functions $\rho^\pm$. Set $t_0=\inf\{t>0: \rho_-(t)>0\}$, so that fibers of $f$ have at most $|B_G(t_0)|$ elements. Let $G/M$ be a generic component of $\square G$, as in Lemma \ref{main}, with $f(G/M)$ contained in a component $H/N$. Then:
\begin{equation}\label{volume2}
\frac{|G/M|}{|B_G(t_0)|}\leq |f(G/M)|\leq |H/N|
\end{equation}
and by formula (\ref{volume}): $|H/N|\leq |B_H(C)|\cdot |G/M|$
This means that, for $n\gg 0$, the set of components with order $n$ of $\square G$ is mapped by $f$ into the set of components of $\square H$ with order between $\frac{n}{|B_G(t_0)|}$ and $|B_H(C)|\cdot n$. 
\hfill$\square$

\section{Examples}

\subsection{Quasi-isometric groups with non-coarsely equivalent full box spaces}

\subsubsection{Connection with the Congruence Subgroup Problem}

{\bf Proof of Theorem \ref{CSP}:} Suppose by contradiction that $\square G$ is coarsely equivalent to $\square H$. By Proposition \ref{subgroupgrowth}, we have $a_n^\triangleleft(G)\leq s_{Bn}^\triangleleft(H)$ for $n\gg 0$.

On the one hand: $a_n^\triangleleft(G)\geq a_n^\triangleleft(\mathbb{F}_2)\geq f(n,2)$ (the number of 2-generated groups of order $n$, up to isomorphism). If $p$ is a prime and $n=p^k$, then $f(p^k,2)\geq p^{ck^2}$ where $k\gg 0$ and $c>0$ is some constant (see Theorem 3.7 in \cite{LubSeg}).

On the other hand: using the CSP for $H$ together with Theorem 6.3.(ii) in \cite{LubSeg}, there exists a constant $K>0$ such that 
$$\log s_{Bn}^\triangleleft(H)\leq K.(\frac{\log Bn}{\log \log Bn})^2.$$

So we get, for $k\gg 0$:
$$c.\log p.k^2\leq K(\frac{\log B+ k.\log p}{\log( \log B+k.\log p)})^2,$$
which is in contradiction with Proposition \ref{subgroupgrowth} for $k\rightarrow\infty$. 
\hfill$\square$

\medskip
In view of the Lubotzky-Zuk conjecture \cite{LubZuk}, $H$ in the above conjecture has CSP if and only if it has property $(\tau)$. In that case, $\square H$ is an expander, while $\square G$ is not, which is another way of proving that $\square G$ and $\square H$ are not coarsely equivalent.

\medskip
Here is an explicit example suggested by N. Nikolov. Let $P(X)\in\Q[X]$ be an irreducible polynomial of degree 3, with one positive root and 2 negative roots (e.g. $P(X)=X^3+11X^2+11X-11$). Let $\alpha$ be the positive root, and let $Q$ be the quadratic form in 3 variables $Q(X_1,X_2,X_3)=X_1^2+X_2^2+\alpha X_3^2$, defined on $\Q(\alpha)$. Let $\mathbb{O}$ be the ring of integers of $\Q(\alpha)$, set $H=SO(Q)(\mathbb{O})$. Then Proposition \ref{CSP} applies to $H$. It is standard to realize $H$ as a co-compact lattice in $SO(2,1)\times SO(2,1)$ (recall the argument: let $\beta_1,\beta_2$ be the negative roots of $P$, and $\sigma_i:\Q(\alpha)\rightarrow \Q(\beta_i)\;(i=1,2)$ be the unique isomorphism mapping $\alpha$ to $\beta_i$; then $(\sigma_1,\sigma_2)$ induces an embedding of $H$ as a lattice in the real points of $SO(\sigma_1(Q))\times SO(\sigma_2(Q))$, which as a Lie group is isomorphic to $SO(2,1)\times SO(2,1)$; co-compactness follows for example from Corollary 5.43 in \cite{Witte}).

Let now $G$ be the product of two surface groups of genus at least 2, viewed as a co-compact lattice in $SO(2,1)\times SO(2,1)$. Then $G$ and $H$ are quasi-isometric, as co-compact lattices in the same Lie groups. Since $G$ maps onto $\mathbb{F}_2$, the box spaces $\square G$ and $\square H$ are not coarsely equivalent, by Theorem \ref{CSP}. Of course this also follows from the fact that $\square H$ is an expander, while $\square G$ is not. 

\subsubsection{Commensurable groups with non-coarsely equivalent full box spaces}

This is based on an idea by K. Bou-Rabee. Let $F$ be the subgroup of $GL_2(\Z)$ generated by $\left(\begin{array}{cc}1 & 0 \\0 & -1\end{array}\right)$ and $\left(\begin{array}{cc}0 & 1 \\1 & 0\end{array}\right)$, so that $F\simeq D_4$ is the isometry group of the square $[-1,1]\times [-1,1]$. Set $G=\Z^2\rtimes F$, and $H=\Z^2$, viewed as an index 8 subgroup in $G$. It is well-known that $a_n^\triangleleft(H)$ is the sum of divisors of $n$ (see e.g. Theorem 15.1 in \cite{LubSeg})

\begin{Lem} $\sqrt{n}\leq s_n^\triangleleft(G)\leq 10\sqrt{n}$.
\end{Lem}

{\bf Proof:} Let $N$ be a finite index subgroup of $G$. Let $p:G\rightarrow F$ be the quotient map. The lattice of normal subgroups of $F$ is very simple:
$$\{Id\}\subset Z(F)\subset F^+\subset F$$
where $|Z(F)|=2$ and $F^+=F\cap SL_2(\Z)$ is cyclic of order 4. So $p(N)$ is one of these 4 subgroups. We look now at $N\cap H$, which we view as an $F$-invariant sub-lattice of $\Z^2$. 

\medskip
{\it Claim:} There exists an integer $k\geq 1$ such that $k\Z^2$ has index 1 or 2 in $N\cap H$ (and $k$ is even in the second case).

To prove the claim, consider the closed cone $C=\{(x,y)\in\R^2: 0\leq y\leq x\}$ in $\R^2$: it is a fundamental domain for the action of $F$ on $\R^2$. Let $v=(x,y)$ be a minimal vector of $N\cap H$, we may assume that $v\in C$. We claim that either $y=0$ or $x=y$. Suppose by contradiction that $0<y<x$. If $y<\frac{\sqrt{3}}{3}.x$, we observe that the length of $v-\left(\begin{array}{cc}1 & 0 \\0 & -1\end{array}\right).v$ (which belongs to $N\cap H$) is shorter than the length of $v$, contradicting minimality; if $y>\frac{\sqrt{3}}{3}.x$, then the length of $v-\left(\begin{array}{cc}0 & 1 \\1 & 0\end{array}\right).v$ (again in $N\cap H$) is shorter than the length of $v$, again a contradiction. 

If $v=(k,0)$ for some $k\geq 1$, then the minimal vectors are $(\pm k,0),(0,\pm k)$ and $N\cap H=k\Z^2$. If $v=(k,k)$ for some $k\geq 1$, then the minimal vectors are the $(\pm k,\pm k)$'s, and $N\cap H$ contains $2k\Z^2$ with index 2.

With this we may count normal subgroups $N$ of given index in $G$: 

\begin{itemize}
\item If $N\cap H=k\Z^2$, the possible values for the index of $N$ in $G$ are: $k^2,\,2k^2,\,4k^2,\,8k^2$ and there is a unique $N$ for each value of the index.
\item If $[N\cap H:2k\Z^2]=2$, the possible values for the index of $N$ in $G$ are $2k^2,\,4k^2,\,8k^2,\,16k^2$ and there is a unique such $N$ for each value of the index. 
\end{itemize}

The result follows.\hfill$\square$

\begin{Cor}\label{noncoarse} In the above example: $\square H$ and $\square G$ are not coarsely equivalent. 
\end{Cor}

{\bf Proof:} Since $a_n^\triangleleft(H)$ grows at least linearly in $n$, and $s_{Bn}^\triangleleft(G)-s_{An}^\triangleleft(G)$ grows at most like $C.\sqrt{n}$, the result follows from Proposition \ref{subgroupgrowth}. \hfill$\square$

\subsection{Non-commensurable groups with coarsely equivalent box spaces}

To construct this example, we will use the following result of Eskin, Fisher and Whyte \cite{EFW13}: let $G$ and $F$ be finite groups; the lamplighter groups $G\wr \mathbb{Z}$ and $F\wr \mathbb{Z}$ are quasi-isometric if and only if there exist positive integers $d,s,r$ such that $|G|=d^s$ and $|F|=d^r$. In fact, we will use the easy direction of the above result, which tells us that $G\wr \mathbb{Z}$ and $F\wr \mathbb{Z}$ admit generating sets with respect to which they are actually isometric if $|G|=|F|$. This is easily achieved by considering a generating set for $G\wr \mathbb{Z}$ consisting of a generator for $\mathbb{Z}$, and the elements of $\bigoplus_{\mathbb{Z}} G$ with an element of $G$ in the 0 position and the identity of $G$ elsewhere, for all elements in $G$. Taking such a generating set for $G\wr \mathbb{Z}$ and $F\wr \mathbb{Z}$, it can be seen that the resulting Cayley graphs are isometric.
 
\begin{Prop} The wreath products $(\mathbb{Z}/4\Z) \wr \mathbb{Z}$ and $(\mathbb{Z}/2\Z \times\mathbb{Z}/2\Z)\wr \mathbb{Z}$ are not commensurable, but admit coarsely equivalent (actually isometric) box spaces.
\end{Prop} 

{\bf Proof:} We first show that the two groups are not commensurable. Given a finite-index subgroup $K$ of the group $(\mathbb{Z}/2\Z\times\mathbb{Z}/2\Z)\wr \mathbb{Z}$, the only elements of finite order contained in $K$ will be elements of order 2. Now consider a finite-index subgroup $H$ of $(\mathbb{Z}/4\Z) \wr \mathbb{Z}$. It necessarily contains elements of order 4, since the index $[\bigoplus_{\mathbb{Z}} \mathbb{Z}/4\Z: (\bigoplus_{\mathbb{Z}} \mathbb{Z}/4\Z)\cap H]$ is finite (since $H$ has finite index in $(\mathbb{Z}/4\Z) \wr \mathbb{Z}$) and so $(\bigoplus_{\mathbb{Z}} \mathbb{Z}/4\Z)\cap H$ cannot be contained in the 2-torsion of $\bigoplus_{\mathbb{Z}} \mathbb{Z}/4\Z$, since the 2-torsion has infinite index in $\bigoplus_{\mathbb{Z}} \mathbb{Z}/4\Z$. Thus $H$ and $K$ cannot be isomorphic, and thus $(\mathbb{Z}/4\Z) \wr \mathbb{Z}$ and $(\mathbb{Z}/2\Z\times\mathbb{Z}/2\Z)\wr \mathbb{Z}$ cannot be commensurable.

We now consider the box spaces of these two groups. Note that they are residually finite, since $\mathbb{Z}$ is residually finite, and the groups $\mathbb{Z}/4\Z$ and $\mathbb{Z}/2\Z\times\mathbb{Z}/2\Z$ are abelian. We will take the box spaces formed of the quotients $(\mathbb{Z}/4\Z) \wr \mathbb{Z}_{2^k}$ and $(\mathbb{Z}/2\Z\times\mathbb{Z}/2\Z)\wr \mathbb{Z}_{2^k}$, for all $k\in \mathbb{N}$. Taking the generating sets of cardinality 5 for the groups $(\mathbb{Z}/4\Z )\wr \mathbb{Z}$ and $(\mathbb{Z}/2\Z\times\mathbb{Z}/2\Z)\wr \mathbb{Z}$ as described above, we see that the box spaces are in fact isometric with respect to the metrics induced by these generating sets. 
\hfill$\square$

\section{Open Problems}
We finish with some questions of interest.

\begin{itemize}
\item
Are the full box spaces $\Box \F_n$, $\Box \F_m$ ($n \neq m$) of free groups of different rank coarsely equivalent?
\item
Given a property $\mathcal{P}$ of box spaces, suppose that $\Box_{(N_k)}G$ has $\mathcal{P}$ with respect to every filtration $(N_k)_{k>0}$ of $G$. When does the full box space $\Box G$ have $\mathcal{P}$?
\item
Consider the wreath product $\Z\wr\Z=\Z[X,X^{-1}]\rtimes\Z$, where $\Z$ acts by multiplication by powers of $X$. Embedding it in a standard way in $GL_2(\Z[X,X^{-1}])$ allows to define congruence subgroups. Is there any congruence box space with property $(D_\alpha)$? Without property $(D_\alpha)$?

\end{itemize}

\medskip
\noindent
Authors' addresses:\\
Institut de Math\'ematiques\\
Unimail\\
11 Rue Emile Argand\\
CH-2000 Neuch\^atel - SWITZERLAND

\medskip
\noindent
anastasia.khukhro@unine.ch\\
alain.valette@unine.ch

\end{document}